\documentclass[11pt]{article}

\usepackage{amssymb,latexsym}
\usepackage{epsfig}
\usepackage{eufrak}
\usepackage{amsmath}
\usepackage{mathrsfs}
\usepackage{color}

\newtheorem{theorem}{Theorem}
\newtheorem{lemma}{Lemma}

\newcommand{\be}{\begin{equation}}
\newcommand{\ee}{\end{equation}}
\newcommand{\bee}{\begin{eqnarray*}}
\newcommand{\eee}{\end{eqnarray*}}
\newcommand{\bel}{\begin{eqnarray}}
\newcommand{\eel}{\end{eqnarray}}
\newcommand{\bec}{\begin{cases}}
\newcommand{\eec}{\end{cases}}
\newcommand{\bem}{\begin{bmatrix}}
\newcommand{\eem}{\end{bmatrix}}

\newcommand{\la}{\label}
\newcommand{\li}{\left}
\newcommand{\ri}{\right}

\newcommand{\vep}{\varepsilon}

\newcommand{\de}{\delta}

\newcommand{\se}{\theta}

\newcommand{\ze}{\zeta}

\newcommand{\f}{\frac}

\newcommand{\cd}{\cdots}

\newcommand{\mscr}{\mathscr}

\newcommand{\mbf}{\mathbf}
\newcommand{\bb}{\mathbb}

\newcommand{\wh}{\widehat}

\newcommand{\mrm}{\mathrm}
\newcommand{\bs}{\boldsymbol}

\newcommand{\tx}{\text}

\newcommand{\iy}{\infty}

\newcommand{\bed}{\begin{description}}
\newcommand{\eed}{\end{description}}
\newcommand{\bei}{\begin{itemize}}
\newcommand{\eei}{\end{itemize}}
\newcommand{\ben}{\begin{enumerate}}
\newcommand{\een}{\end{enumerate}}
\newcommand{\bib}{\bibitem}
\newcommand{\beL}{\begin{lemma}}
\newcommand{\eeL}{\end{lemma}}
\newcommand{\beT}{\begin{theorem}}
\newcommand{\eeT}{\end{theorem}}

\newcommand{\bpf}{\begin{pf}}
\newcommand{\epf}{\end{pf}}
\newcommand{\bsk}{\bigskip}

\newcommand{\pfbox}{\hfill\mbox{$\Box$}}

\newenvironment{pf}{\paragraph*{Proof{\rm.}}}{\pfbox\bigskip}

\begin{document}

\title{{\bf Confidence Interval for the
Mean of a Bounded Random Variable and Its Applications in Point
Estimation}
\thanks{The author is currently with Department of Electrical Engineering,
Louisiana State University at Baton Rouge, LA 70803, USA, and Department of Electrical Engineering, Southern University and A\&M College, Baton
Rouge, LA 70813, USA; Email: chenxinjia@gmail.com}}

\author{Xinjia Chen}

\date{November, 2010}

\maketitle

\begin{abstract}

In this article, we derive an explicit formula for computing
confidence interval for the mean of a bounded random variable.
Moreover, we have developed multistage point estimation methods for
estimating the mean value with prescribed precision and confidence
level based on the proposed confidence interval.

\end{abstract}

\section{Introduction}

In many areas of sciences and engineering, it is a frequent problem
to estimate  the mean of a bounded random variable. Conventional
technique for constructing confidence interval relies on the Central
Limit Theorem. However, for small and moderate sample size, using
normal approximation can lead to serious under-coverage of the mean.
In the case of bounded random variables, even the sample size is
very large, the error can also be intolerable when the parent
distribution is highly skewed toward extremes.

In this article, by applying an inequality obtained by Massart 1990
and Hoeffding's probability inequality, we have derived an explicit
formula for interval estimation of the mean in the bounded case. The
formula is extremely simple.  Moreover, we have proposed multistage
estimation methods for estimating the mean value with prescribed
precision and confidence level based on the construction of
confidence interval.

\section{Explicit Formula}

Since any random variable $X$ bounded in interval $[a, b]$ (i.e., $\Pr \{a \leq X \leq b \} = 1$)
 has a linear relation with random variable $Z = \frac{X-a}{b-a}$, it suffices to consider interval estimation for the mean of
 random variable $Z$ on interval $[0,1]$ (i.e., $\Pr \{0 \leq Z \leq 1 \} = 1$) and employ transformation
 $X = (b-a) Z + a$ to obtain an estimation for the mean of $X$.
 The following Theorem 1 provides an easy method for constructing confidence interval for the mean of $Z$.

\begin{theorem} \label{continuous}
Let $\delta \in (0,1)$ and $c = \frac{9}{ 2 \ln \frac{2}{\delta} }$.
Let $\Pr \{0 \leq Z \leq 1 \} = 1$ and $\mu = \bb{E}(Z)$. Let
$\overline{Z} = \frac{\sum_{i=1}^n Z_i}{n}$ where $n$ is the sample
size and $Z_i, \;\; i = 1, \cdots, n$ are i.i.d. observations of
$Z$.  Define
\[
L = \overline{Z} + \frac{3}{4 + n c} \left[ 1 - 2 \overline{Z} -
\sqrt{ 1 + n c \overline{Z} (1 - \overline{Z}) } \right],
\]
\[
U = \overline{Z} + \frac{3}{4 + n c} \left[ 1 - 2 \overline{Z} +
 \sqrt{ 1 + n c \overline{Z} (1 - \overline{Z}) } \right].
\]
Then,
\[
\Pr \{ L < \mu < U \} \geq 1 - \delta.
\]
\end{theorem}

To prove Theorem 1, we need some preliminary lemmas.

\begin{lemma}
Let $\alpha = \frac{1}{n c}$.  Let $0 \leq t \leq 1$. Then
$\epsilon(t) = \frac{ 3 \alpha ( 1 - 2 t) + 3 \sqrt{\alpha^2 + 4
\alpha t (1 - t) } } { 2(1 + \alpha)} \geq 0$ satisfies equation \be
\exp \left ( - \frac{ n \epsilon^2 } { 2 ( t + \frac{\epsilon} {3})
( 1 -  t - \frac{\epsilon} {3} ) } \right ) = \frac{\delta}{2}
\label{eq1} \ee with respect to $\epsilon$.
\end{lemma}

\begin{pf}
Let $q = t + \frac{\epsilon}{3}$ where $\epsilon$ satisfies equation~(\ref{eq1}).
Then $q$ satisfies equation
$
\exp \left ( - \frac{ 9 n (q-t)^2 }
{ 2 q
( 1 -  q ) } \right ) = \frac{\delta}{2}$,
which can be simplified as
\be
(q - t)^2 + \alpha q (q -1) = 0
\label{eq2}
\ee
with two real roots
$q = \frac{ 2 t + \alpha \pm \sqrt{\alpha^2 + 4 \alpha t (1 - t) } } { 2(1 + \alpha)}$.
Making use of the relation between $\epsilon$ and $q$, we find the roots of equation (\ref{eq1}) as
$\epsilon_1 = \frac{ 3 \alpha ( 1 - 2 t) + 3 \sqrt{\alpha^2 +
4 \alpha t (1 - t) } } { 2(1 + \alpha)}$
and $\epsilon_2 = \frac{ 3 \alpha ( 1 - 2 t) - 3 \sqrt{\alpha^2 +
4 \alpha t (1 - t) } } { 2(1 + \alpha)}$.
It can be verified that
$| \alpha ( 1 - 2 t) | ^2 \leq \alpha^2 +
4 \alpha t (1 - t)$, which leads to $\epsilon(t) = \epsilon_1 \geq 0$ and $\epsilon_2 \leq 0$.
\end{pf}

\begin{lemma} Let $t \in (0,1)$.
Then $\epsilon(t)$
is a concave function with respect to $t$.
\end{lemma}

\begin{pf}
By equation (\ref{eq2}), we have
$0 < t < q < 1$ and
$
\frac{d q  } {d t } = \frac{ 2 (q-t)  } {2 (q-t) + \alpha ( 2q -1) }
= \frac{ 1   } {  1 + \alpha + \frac{ \alpha (t- \frac{1}{2}) } {q - t} }$.  Consequently,
$
\frac{d \left( \frac{t-\frac{1}{2}}{q-t}  \right)   } { dt  } > 0 \;\;\Longleftrightarrow \;\;
(q - t) - (t - \frac{1} {2} ) (  \frac{d q  } {d t } - 1  ) > 0
\;\;\Longleftrightarrow \;\;
q - t > \frac{ \alpha (t- \frac{1}{2}) (1 - 2 q)} { 2 (q-t) + \alpha ( 2q -1) }.
$
Moreover,
$
\frac{d^2 \epsilon } {d t^2 } = 3 \frac{d^2 q  } {d t^2 }
= \frac{ -3\alpha   } {  \left[1 + \alpha +
\frac{ \alpha (t- \frac{1}{2}) } {q - t} \right]^2} \;
 \frac{d \left( \frac{t-\frac{1}{2}}{q-t}  \right)   } { dt  }$.
 Therefore, to show $\frac{d^2 \epsilon(t) } {  d t^2  } < 0$,
 it suffices to show inequality $q - t > \frac{ \alpha
 (t- \frac{1}{2}) (1 - 2 q)} { 2 (q-t) + \alpha ( 2q -1) }$, which is equivalent to
$1 > \frac{ \alpha (t- \frac{1}{2}) (1 - 2 q)}
{ 2 (q-t)^2 + \alpha (q-t) ( 2q -1) }$ since $q-t > 0$.
Note that
$\frac{ \alpha (t- \frac{1}{2}) (1 - 2 q)}
{ 2 (q-t)^2 + \alpha (q-t) ( 2q -1) } = \frac{ \alpha (t- \frac{1}{2}) (1 - 2 q)}
{ 2 (q-t)^2 + 2 \alpha q (q-1)  + \alpha q - \alpha t ( 2q -1) }
= \frac{ (t- \frac{1}{2}) (1 - 2 q)}
{ q -  t ( 2q -1) }$
because $q$ satisfies equation (\ref{eq2}).  It follows that, to show $\frac{d^2 \epsilon(t) } {  d t^2  } < 0$,
 it suffices to show inequality $1 > \frac{ (t- \frac{1}{2}) (1 - 2 q)} { q -  t ( 2q -1) }$.
Invoking inequality  $0 < t < q < 1$, we can show that $q -  t ( 2q -1) > 0$, which leads to
equivalent relations
$
1 > \frac{ (t- \frac{1}{2}) (1 - 2 q)} { q -  t ( 2q -1) }  \;\;\Longleftrightarrow \;\;
 q -  t ( 2q -1) > (t- \frac{1}{2}) (1 - 2 q)
\;\;\Longleftrightarrow \;\;
0 > - \frac{1}{2}.
$
The last inequality is trivially true.
\end{pf}

\begin{lemma}
Let $\beta = \frac{4}{n c}$.  Let $t(z) = z + \frac{ 3\beta(1-2z) -
 3 \sqrt{\beta^2 + 4 \beta z (1 - z) } } { 4(1 + \beta)}$ where $0 \leq z \leq 1$.
Then $z - t(z) = \epsilon(t(z))$ and $t(z) \leq z$.
\end{lemma}

\begin{pf}

Let $p = t + \frac{z-t}{3}$ where $t$ satisfies $z-t=\epsilon(t)$.
It follows that $\epsilon(t) = \frac{-3(p-z)}{2}$ and $t + \frac{\epsilon(t)} { 3} = p$.
By Lemma 1, $\epsilon(t)$ satisfies equation (\ref{eq1}), hence $p$ satisfies equation
$\exp \left ( - \frac{ \frac{9}{4} n (p-z)^2 }
{ 2 p( 1 -  p) } \right ) = \frac{\delta}{2}$, which can be simplified as $(p - z)^2 + \beta p (p-1) = 0$ with
two roots
$p = \frac{ 2 z + \beta \pm \sqrt{\beta^2 + 4 \beta z (1 - z) } } { 2(1 + \beta)}$.
Making use of the relation between $p$ and $t$, we find the solution of equation
$z-t=\epsilon(t)$ with respect to $t$
as $t_1 = z + \frac{ 3\beta(1-2z) + 3 \sqrt{\beta^2 + 4 \beta z (1 - z) } } { 4(1 + \beta)}$ and
$t_2 = z + \frac{ 3\beta(1-2z) - 3 \sqrt{\beta^2 + 4 \beta z (1 - z) } } { 4(1 + \beta)}$.  It can be shown that
$|\beta(1-2z)|^2 \leq \beta^2 + 4 \beta z (1 - z)$, which leads to $t_1 \geq z$ and $t_2 \leq z$.
So the proof is completed by noting that $t(z) = t_2$.
\end{pf}

\begin{lemma} Let $0 < \mu < 1$ and $0 \leq z \leq 1$.  Then $z - \mu \geq \epsilon(\mu)$ if $t(z) \geq \mu$.
\end{lemma}

\begin{pf}
Let $t(z) \geq \mu > 0$.  By Lemma 3, we have $z - t(z) \geq 0$ and thus $z - \mu \geq z - t(z) \geq 0$.
We claim that $z - \mu > 0$. If this is not true, then $z = \mu$ and $t(z) \geq z > 0$.
By Lemma 3, we have $t(z) = z > 0$. On the other hand, $t(z) = z$ results in $z = 0$.
Thus we arrive at contradiction $0 > 0$.
So we have shown $z - \mu > 0$ and it follows that
$0 \leq \frac{z - t(z) } {z - \mu } \leq 1$.  We next show that $z - \mu \geq \epsilon(\mu)$.
Suppose for the purpose of contradiction that $z - \mu < \epsilon(\mu)$.
Then
\[
z - t(z) = (z - \mu ) \frac{z - t(z) } {z - \mu } <
\epsilon( \mu )  \frac{z - t(z) } {z - \mu } + \left( 1 - \frac{z - t(z) } {z - \mu } \right) \epsilon(z).
\]
By Lemma 2, $\epsilon(t)$ is concave with respect to $t$, hence
$\epsilon( \mu ) \frac{z - t(z) } {z - \mu }  + ( 1 - \frac{z - t(z) } {z - \mu } ) \epsilon(z) < \epsilon(t(z))$, which yields
$z - t(z) < \epsilon(t(z))$.  Recall Lemma 3, $z - t(z) = \epsilon(t(z))$.
It follows that $\epsilon(t(z)) < \epsilon(t(z))$, which is a contradiction.

\end{pf}

We are now in the position to prove Theorem 1.  By Theorem 1 of Hoeffding 1963, \be \Pr \{ \overline{Z} \geq \mu
+ \epsilon \} \leq \left\{ \li (  \frac{\mu} {\mu + \epsilon}  \ri )^{ \mu + \epsilon } \; \li (  \frac{1- \mu}
{1- \mu - \epsilon}  \ri )^{ 1- \mu - \epsilon } \right\}^n \;\;\;\; \forall \epsilon \in (0, 1 - \mu).
\label{h} \ee By Lemma 1 of  Massart 1990, \be (\mu + \epsilon) \ln \li (\frac{\mu+\epsilon}{\mu} \ri ) + (1 -
\mu - \epsilon) \ln \li (\frac{1 - \mu - \epsilon}{1 - \mu} \ri )
 \geq \frac{ \epsilon^2 }
 {2 (\mu + \frac{\epsilon} { 3 }) (1 - \mu - \frac{\epsilon} { 3 })  } \;\;\;\; \forall \epsilon \in (0, 1 - \mu).
\label{m} \ee It follows from (\ref{h}) and (\ref{m}) that \be \Pr
\{  \overline{Z} \geq \mu + \epsilon  \} \leq \exp \left ( - \frac{
n \epsilon^2 } { 2 ( \mu + \frac{\epsilon} {3}) ( 1 -  \mu -
\frac{\epsilon} {3} ) } \right )  \;\;\;\; \forall \epsilon > 0.
\label{mh} \ee By the definition of $t(.)$, we can verify that $L =
t( \overline{Z} )$.  Thus $\Pr \{ L \geq \mu   \} = \Pr \{  t(
\overline{Z} ) \geq \mu  \}$. Applying Lemma 4, we have $\Pr \{  t(
\overline{Z} ) \geq \mu  \} \leq \Pr \{  \overline{Z} - \mu \geq
\epsilon( \mu)  \}$.  Hence by (\ref{mh}) and Lemma 1,
\[
\Pr \{ L \geq \mu   \} \leq \Pr \{  \overline{Z} - \mu \geq
\epsilon( \mu)  \} \leq \exp \left ( - \frac{ n [\epsilon( \mu )]^2
} { 2 ( \mu + \frac{\epsilon ( \mu )} {3}) ( 1 -  \mu -
\frac{\epsilon ( \mu )} {3} ) } \right ) = \frac{\delta}{2}.
\]
Since $\Pr \{ L \geq \mu   \} \leq \frac{\delta}{2}$ has been shown,
applying this conclusion to random variable $1 - Z$, we have  $ \Pr
\{ U \leq \mu   \} \leq \frac{\delta}{2}. $

Finally, by applying Bonferrnoni's inequality,  we have
\begin{eqnarray*}
\Pr \{ L < \mu < U  \} & \geq & \Pr \{ L < \mu   \} + \Pr \{ U > \mu   \} -1\\
& = & 1 - \Pr \{ L \geq \mu   \} + 1 - \Pr \{ U \leq \mu   \}  -1\\
& \geq & 1 - \frac{\delta}{2} + 1 - \frac{\delta}{2} -1 = 1 -\delta.
\end{eqnarray*}

\section{Applications in Multistage Point Estimation}

We would like to point out that the simple interval estimation
method described above can be used to construct multistage sampling
plans for estimating the mean value of a bounded variable with
prescribed precision and confidence level.   To illustrate such
applications, we shall first present some general results of
multistage point estimation based on confidence intervals.

Let $X$ be a random variable parameterized by $\se$, which is not
necessary bounded. Let $X_1, X_2, \cd$ be a sequence of random
samples of $X$. The goal is to estimate $\se$ via a multistage
sampling plan with the following structure. The sampling process is
divided into $s$ stages, where $s$ can be infinity or a positive
integer. The continuation or termination of sampling is determined
by decision variables.  For each stage with index $\ell$, a decision
variable $\bs{D}_\ell = \mscr{D}_\ell (X_1, \cd, X_{\mbf{n}_\ell})$
is defined based on samples $X_1, \cd, X_{\mbf{n}_\ell}$, where
$\mbf{n}_\ell$ is the number of samples available at the $\ell$-th
stage.  It should be noted that $\mbf{n}_\ell$ can be a random
number, depending on specific sampling schemes. The decision
variable $\bs{D}_\ell$ assumes only two possible values $0, \; 1$
with the notion that the sampling is continued until $\bs{D}_\ell =
1$ for some $\ell$.  For the $\ell$-th stage, an estimator
$\wh{\bs{\se}}_\ell$ for $\se$ is defined based on samples $X_1,
\cd, X_{\mbf{n}_\ell}$. Let $\bs{l}$ denote the index of stage when
the sampling is terminated.  Then, the point estimator for $\se$,
denoted by $\wh{\bs{\se}}$, is equal to  $\wh{\bs{\se}}_{\bs{l}}$.
The decision variables $\bs{D}_\ell$ can be defined in terms of
estimators $\wh{\bs{\se}}_\ell$ and confidence intervals $(L_\ell,
U_\ell)$, where the lower confidence limit $L_\ell$ and upper
confidence limit $U_\ell$ are functions of $X_1, \cd,
X_{\mbf{n}_\ell}$ for $\ell = 1, \cd, s$.  Depending on various
error criterion, we have different sampling plans as follows.

\beT \la{thm2} Let $\vep > 0, \; \ze > 0$ and $\de \in (0, 1)$.  For
$\ell = 1, \cd, s$, let $(L_\ell, U_\ell)$ be a confidence interval
such that $\Pr \{ L_\ell < \se < U_\ell \} > 1 - \ze \de$.   Suppose
the stopping rule is that sampling is continued until $U_\ell - \vep
< \wh{\bs{\se}}_\ell < L_\ell + \vep$ at some stage with index
$\ell$. Then, $\Pr \{ | \wh{\bs{\se}} - \se | < \vep \} > 1 - \de$
provided that $ s \ze < 1$ and that $\Pr \{ U_s - \vep <
\wh{\bs{\se}}_s < L_s + \vep \} = 1$. \eeT

\bpf

By the assumption that $\Pr \{ U_s - \vep < \wh{\bs{\se}}_s < L_s +
\vep \} = 1$, we have that $\Pr \{ \bs{l} > s \} = 0$.  Hence, by
the definition of the sampling scheme described by Theorem 2, we
have {\small \bee \Pr \{ | \wh{\bs{\se}} - \se | \geq \vep \} & = &
\sum_{\ell = 1}^s \Pr \{ | \wh{\bs{\se}}_\ell - \se | \geq \vep, \;
\bs{l} = \ell \} \leq \sum_{\ell = 1}^s \Pr \{ | \wh{\bs{\se}}_\ell
- \se |
\geq \vep, \; \bs{D}_\ell = 1 \}\\
& \leq  & \sum_{\ell = 1}^s \Pr \{  | \wh{\bs{\se}}_\ell - \se |
\geq \vep, \; U_\ell - \vep < \wh{\bs{\se}}_\ell < L_\ell + \vep
\}\\
& =  & \sum_{\ell = 1}^s \Pr \{  \wh{\bs{\se}}_\ell \geq \se + \vep
\; \; \tx{or} \; \; \wh{\bs{\se}}_\ell \leq \se - \vep,
\; U_\ell - \vep < \wh{\bs{\se}}_\ell < L_\ell + \vep \} \\
& \leq  & \sum_{\ell = 1}^s \Pr \{  L_\ell > \wh{\bs{\se}}_\ell -
\vep \geq \se \; \; \tx{or} \; \;
 U_\ell < \wh{\bs{\se}}_\ell + \vep \leq \se  \}\\
& \leq  & \sum_{\ell = 1}^s \Pr \{  L_\ell \geq \se  \; \; \tx{or}
\; \; U_\ell \leq \se \} =   \sum_{\ell = 1}^s  \li [ 1 - \Pr \{
L_\ell < \se < U_\ell \} \ri ].  \eee} Therefore, by the assumption
that $\Pr \{ L_\ell < \se < U_\ell \} > 1 - \ze \de$ for $\ell = 1,
\cd, s$, we have $\Pr \{ | \wh{\bs{\se}} - \se | \geq \vep \} \leq
\sum_{\ell = 1}^s  \li [ 1 - \Pr \{ L_\ell < \se < U_\ell \} \ri ] <
s \ze \de$, from which the theorem immediately follows.

\epf

Theorem 2 indicates that the coverage probability $\Pr \{ |
\wh{\bs{\se}} - \se | < \vep \}$ can be adjusted by $\ze > 0$. In
order to make the coverage probability above $1 - \de$, it suffices
to choose a sufficiently small $\ze > 0$. We would like to point out
that, for estimating the mean value of a random variable bounded in
$[a, b]$, Theorem 2 can be applied based on the following choice:

(i) The sample sizes of the sampling plan are chosen as
deterministic integers $n_1 < \cd < n_s$ such that $n_s > \f{(b -
a)^2} { 2 \vep^2} \ln \f{2} {\ze \de}$.

(ii) The confidence intervals are constructed by virtue of Theorem
1.

\bsk

 Theorem \ref{thm3} at below describes a method for defining
stopping rules  for estimating $\se$ with relative precision so that
the coverage probabilities can be controlled by $\ze$.

\beT \la{thm3} Let $\vep > 0, \; \ze > 0$ and $\de \in (0, 1)$.  For
$\ell = 1, \cd, s$, let $(L_\ell, U_\ell)$ be a confidence interval
such that $\Pr \{ L_\ell < \se < U_\ell \} > 1 - \ze \de$.   Suppose
the stopping rule is that sampling is continued until $[1 -
\mrm{sgn} ( \wh{\bs{\se}}_\ell ) \; \vep ] U_\ell <
\wh{\bs{\se}}_\ell < [1 + \mrm{sgn} ( \wh{\bs{\se}}_\ell ) \; \vep ]
L_\ell$ at some stage with index $\ell$. Then, $\Pr \{ |
\wh{\bs{\se}} - \se | < \vep |\se| \} > 1 - \de$ provided that $ s
\ze < 1$ and that $\Pr \{ [ 1 - \mrm{sgn} ( \wh{\bs{\se}}_s ) \;
\vep ] U_s < \wh{\bs{\se}}_s < [ 1 + \mrm{sgn} ( \wh{\bs{\se}}_s )
\;  \vep ] L_s \} = 1$, where $\mrm{sgn} (x)$ is the sign function
which assumes values $1, \; 0$ and $-1$ for $x > 0, \; x = 0$ and $x
< 0$ respectively.  \eeT

We would like to note that, for estimating the mean value of a
random variable bounded in $[0, 1]$, we can use Theorems 1 and 3
based on multistage inverse sampling.

\bsk

Theorem \ref{thm4} at below describes a method for defining stopping
rules  for estimating $\se$ with mixed precision so that the
coverage probabilities can be controlled by $\ze$.

\beT \la{thm4}  Let $0 < \de < 1, \; \vep_a
> 0, \; \vep_r > 0$ and $\ze > 0$.  For $\ell = 1, \cd, s$, let $(L_\ell, U_\ell)$
be a confidence interval such that $\Pr \{ L_\ell < \se < U_\ell \}
> 1 - \ze \de$.  Suppose the stopping rule is that sampling
is continued until $U_\ell - \max ( \vep_a, \; \mrm{sgn}
(\wh{\bs{\se}}_\ell) \; \vep_r U_\ell ) < \bs{\wh{\se}}_\ell <
L_\ell + \max ( \vep_a, \; \mrm{sgn} (\wh{\bs{\se}}_\ell) \; \vep_r
L_\ell )$ at some stage with index $\ell$.
 Then, {\small $\Pr \li \{ \li | \bs{\wh{\se}} - \se \ri | < \vep_a \; \tx{or} \;
  \li | \bs{\wh{\se}} - \se \ri | < \vep_r |\se| \ri \}
\geq 1 - \de$} provided that $ s \ze < 1$ and that $\Pr \{ U_s -
\max ( \vep_a, \; \mrm{sgn} (\wh{\bs{\se}}_s) \; \vep_r U_s ) <
\bs{\wh{\se}}_s < L_s + \max ( \vep_a, \; \mrm{sgn}
(\wh{\bs{\se}}_s) \; \vep_r L_s )  \} = 1$. \eeT

For estimating the mean value of a random variable bounded in $[a,
b]$, Theorem 4 can be applied based on the following choice:

(i) The sample sizes of the sampling plan are chosen as
deterministic integers $n_1 < \cd < n_s$ such that $n_s > \f{(b -
a)^2} { 2 \vep^2} \ln \f{2} {\ze \de}$.

(ii) The confidence intervals are constructed by virtue of Theorem
1.

\bsk

In Theorems 2--4, the number of stages, $s$, is assumed to be a
finite integer.  In some situations, a sampling plan with a finite
number of stages is impossible to guarantee the prescribed precision
and confidence level.  In this regard, the following theorems are
useful.

\beT \la{thm5} Let $\vep > 0, \; \ze > 0$ and $\de \in (0, 1)$.  Let
$\tau$ be a positive integer.  Let $(L_\ell, U_\ell)$ be a
confidence interval such that $\Pr \{ L_\ell < \se < U_\ell \} > 1 -
\ze \de$ for $\ell \leq \tau$ and that $\Pr \{ L_\ell < \se < U_\ell
\} > 1 - \ze \de 2^{\tau - \ell}$ for $\ell > \tau$. Suppose the
stopping rule is that sampling is continued until $U_\ell - \vep <
\wh{\bs{\se}}_\ell < L_\ell + \vep$ at some stage with index $\ell$.
Then, $\Pr \{ | \wh{\bs{\se}} - \se | < \vep \} > 1 - \de$ provided
that $ (\tau + 1) \ze < 1$ and that $\Pr \{ \bs{l} < \iy  \} = 1$.
\eeT

\beT \la{thm6} Let $\vep > 0, \; \ze > 0$ and $\de \in (0, 1)$.  Let
$\tau$ be a positive integer.  Let $(L_\ell, U_\ell)$ be a
confidence interval such that $\Pr \{ L_\ell < \se < U_\ell \} > 1 -
\ze \de$ for $\ell \leq \tau$ and that $\Pr \{ L_\ell < \se < U_\ell
\} > 1 - \ze \de 2^{\tau - \ell}$ for $\ell > \tau$. Suppose the
stopping rule is that sampling is continued until $[1 - \mrm{sgn} (
\wh{\bs{\se}}_\ell ) \; \vep ] U_\ell < \wh{\bs{\se}}_\ell < [1 +
\mrm{sgn} ( \wh{\bs{\se}}_\ell ) \; \vep ] L_\ell$ at some stage
with index $\ell$. Then, $\Pr \{ | \wh{\bs{\se}} - \se | < \vep
|\se| \}
> 1 - \de$ provided that $(\tau + 1) \ze < 1$ and that $\Pr \{
\bs{l} < \iy  \} = 1$. \eeT

\beT \la{thm7} Let $0 < \de < 1, \; \vep_a
> 0, \; \vep_r > 0$ and $\ze > 0$.  Let $\tau$ be a positive integer.  Let
$(L_\ell, U_\ell)$ be a confidence interval such that $\Pr \{ L_\ell
< \se < U_\ell \} > 1 - \ze \de$ for $\ell \leq \tau$ and that $\Pr
\{ L_\ell < \se < U_\ell \} > 1 - \ze \de 2^{\tau - \ell}$ for $\ell
> \tau$.  Suppose the stopping rule is that sampling
is continued until $U_\ell - \max ( \vep_a, \; \mrm{sgn}
(\wh{\bs{\se}}_\ell) \; \vep_r U_\ell ) < \bs{\wh{\se}}_\ell <
L_\ell + \max ( \vep_a, \; \mrm{sgn} (\wh{\bs{\se}}_\ell) \; \vep_r
L_\ell )$ at some stage with index $\ell$.
 Then, {\small $\Pr \li \{ \li | \bs{\wh{\se}} - \se \ri | < \vep_a \; \tx{or} \;
  \li | \bs{\wh{\se}} - \se \ri | < \vep_r |\se| \ri \}
\geq 1 - \de$} provided that $(\tau + 1) \ze < 1$ and that $\Pr \{
\bs{l} < \iy  \} = 1$. \eeT

We would like to note that, for estimating the mean value of a
random variable bounded in $[a, b]$, Theorems 5--7 can be used since
it can be shown that $\Pr \{ \bs{l} < \iy  \} = 1$ as a consequence
of using the confidence interval described by Theorem 1.  In this
paper, we have omitted the proofs of Theorems 3--7, since these
theorems can be readily shown in a more general setting by virtue of
identity (1) and Theorem 3 in the 22th version of our paper
\cite{Chen_EST}.  Although theorems 2--7 propose general methods to
define stopping rules so that the associated coverage probabilities
can be controlled by $\ze$, no specific method is provided for using
$\ze$ to adjust the coverage probabilities as close as possible to
the desired level $1 - \de$.  This issue is extensively explored in
our paper \cite{Chen_EST}.

\end{document}